\documentstyle[12pt]{article}

\newlength{\abstractwidth}
\renewcommand{\thefootnote}{\fnsymbol{footnote}}
\renewcommand{\thanks}[1]{\footnote{#1}} 
\newcommand{\starttext}{ \setcounter{footnote}{0}
\renewcommand{\thefootnote}{\arabic{footnote}}}

\newcommand{\be}{\begin{equation}}
\newcommand{\bea}{\begin{eqnarray}}
\newcommand{\eea}{\end{eqnarray}} \newcommand{\ee}{\end{equation}}
\newcommand{\N}{{\cal N}} 
 \def\ba{\begin{eqnarray}}
\def\ea{\end{eqnarray}}

\def\D{{\cal D}}

\def\N{{\cal N}}

\def\o{\omega}

\def\tr{{\rm tr}}

\def\log{\,{\rm log}\,}

\def\o{\omega}

\def\a{\alpha}

\def\b{\beta}
\def\g{\gamma}
\def\d{\delta}

\def\l{\lambda}

\def\o{\omega}

\def\D{\Delta}
\def\O{\Omega}

\def\vp{\varphi}

\def\na{\nabla}

\def\ve{\varepsilon}

\def\N{\bf N}

\def\R{{\bf R}}

\def\i{\infty}

\def\p{\partial}

\def\D{\Delta}

\def\na{{\nabla}}

\def\[{{\bf [}}
\def\]{{\bf ]}}

{
	\par\vspace{\baselineskip}%
	\textbf{Proof}\begin{itshape}%
		\par\vspace{\baselineskip}%
	}%
	{\end{itshape}}

\begin{document}
\starttext \baselineskip=18pt \setcounter{footnote}{0}
\newtheorem{theorem}{Theorem}
\newtheorem{lemma}{Lemma}
\newtheorem{corollary}{Corollary}
\newtheorem{definition}{Definition}
\newtheorem{conjecture}{Conjecture}
\newtheorem{proposition}{Proposition}
\newtheorem{example}{Example}
\newtheorem{remark}{Remark}

\begin{center}
{\Large \bf Convergence of Type IIA manifolds and application to the Type IIA flow.
}

\medskip
\centerline{Nikita Klemyatin\footnote{Supported in part by the National Science Foundation under grant DMS-1855947.}}

\medskip

\begin{abstract}

{\footnotesize }

The Type IIA flow is a flow of 3-forms on a 6-dimensional symplectic manifold. It was introduced by Fei, Phong, Picard, and Zhang in \cite{FPPZb}. There is little known about the singularities of this flow in general. We formulate and prove convergence theorems for this flow. We also describe the singularity models for this flow.

\end{abstract}

\end{center}

\baselineskip=15pt
\setcounter{equation}{0}
\setcounter{footnote}{0}


\section{Introduction}
\setcounter{equation}{0}

The problem of finding canonical geometric structures has a long history with roots both in geometry and physics. There are many examples of such geometric structures: Einstein and K\"ahler-Einstein metrics on Riemannian and K\"ahler manifolds, Hermitian-Einstein metrics on stable holomorphic vector bundles (see \cite{UY}), and Ricci-flat metrics on Calabi-Yau manifolds (\cite{Y}). However, many more examples of canonical geometric structures arise from the string theory. Starting with the celebrated work of Candelas, Horowitz, Strominger, and Witten \cite{CHSW}, various geometric structures arise from the string theory were studied extensively (see \cite{FGP, _Fei15_,_Fei16_, FY, _FP19_, _FPic19_, _FP21,FPPZa,LY,PPZ16,PPZ17,PPZAnn,B,BX} for some examples). The Type IIA geometry is one such structure. It was introduced by Fei, Phong, Picard, and Zhang in \cite{FPPZb} as a geometric structure emerging from the Type IIA string theory (see \cite{TS1}, \cite{TS2} and \cite{TS3}). This is a specific geometrical structure on a 6-dimensional symplectic manifold $(M,\o)$ together with a primitive 3-form $\vp$ on $M$, which satisfies an additional nondegeneracy condition. Such form $\vp$ defines an almost complex structure $J_\vp$ on $M$, which is compatible with $\o$ and satisfies the following system of equations:

\bea
d\Lambda d(|\vp|^2 \star \vp) = \rho, 
\\ 
d\vp = 0.
\eea
Here $\Lambda$ is the usual contraction with the symplectic form $\o$, $\star$ is the usual Hodge-star operator with respect to the metric $g_\vp$ defined by $\o$ and $J_\vp$. Also, $\rho$ is a source form, which is Poincare dual to a linear combination of special Lagrangian cycles, calibrated by $\vp$.

In \cite{FPPZb} the authors suggest a geometric flow approach to solve this system. Namely, they considered the following flow of closed, primitive 3-forms on a 6-dimensional symplectic manifold $(M, \o)$:

\bea
\p_t \vp = d\Lambda d(|\vp|^2 \star \vp) - \rho.
\eea

One can check that this flow preserves the primitiveness and the closedness of $\vp$ if the initial data is primitive.

In \cite{FPPZb} the authors established some important properties of this flow, including the short-time existence, the uniqueness, and the Shi-type estimates for the source-free flow (i.e. $\rho = 0$). In other papers \cite{FPPZc}, \cite{FPPZd} and \cite{_FP21_} authors established additional properties of the flow, such as dynamical stability on compact symplectic 6-manifolds, as well as applications to symplectic geometry. However, it is little known about the formation of singularities of the Type IIA flow. Once we have the Shi-type estimates, we can study the formation of singularities, as it was done in the case of the Ricci flow by Hamilton \cite{Hamilton1}, for the Laplacian flow in \cite{LotayWei15}, and in the case of the Type IIB flow in \cite{NK}.

In order to do it, we need to establish the compactness properties of the solutions to the Type IIA flow. First, we need to establish convergence of manifolds, which carry the Type IIA manifolds:

\begin{theorem}\label{CompactnessTypeIIA}
	Let $\{(M_j, g_j, \o_j, \vp_j, J_j, p_j)\}^\infty_{j=1}$ be a sequence of complete, pointed Type IIA manifolds. Denote $u_j := \log |\phi_j|_{g_j}$. Assume there are positive constants $C_k, \iota_0 > 0$, independent of $j$, such that the following holds:
	
	\smallskip
	
	{\rm(1)} 
	\bea
	|(\na^{g_j})^k Rm_j|_{g_j} \leq C_k, 
	\\
	|(\na^{g_j})^k N_j|_{g_j} \leq C_k, 
	\\
	|(\na^{g_j})^k \vp_j|_{g_j} \leq C_k,
	\\
	|(\na^{g_j})^k u_j|_{g_j} \leq C_k
	\eea 
	
	{\rm(2)} The injectivity radius satisfies the following inequality:
	\bea
	\mathrm{inj}(M_j, g_j) \geq \iota_0
	\eea
	
	\smallskip
	
	Then there is a complete, pointed manifold $(M, g, \o, \vp, J, p)$,  a subsequence \\ $\{(M_{j_k}, g_{j_k}, \o_{j_k}, \vp_{j_k}, J_{j_k}, p_{j_k})\}$ which converge to $(M, g, \omega, \phi, J, p)$ in the sense of Definition \ref{GCH_Convergence} below.
\end{theorem}

Once we have Theorem \ref{CompactnessTypeIIA}, we can prove the theorem about the convergence of solutions to the Type IIA flow:

\begin{theorem}\label{CompactnessTypeIIAFlow}
	Let $(M_j,g_j(t),\o_j,\vp_j(t),J_j(t),p_j)$ be the family of complete solutions of the Type IIA flow for $t \in (\a;\b) \subset \R$. Assume that the following holds:
	
	\smallskip
	
	{\rm(1)} There is a constant $\iota_0>0$ such that $\mathrm{inj}(M_j, g_j(t)) \geq \iota_0$ along the flow;
	
	{\rm(2)}
	There is a constant $C>0$ indepentent of $j$, such that $|u_j(t)| + |Rm_j(t)| \leq C $.
	
	\smallskip
	
	Then there is a complete solution $(M, g(t), \o, \vp(t), J(t), p)$ and a subsequence	\\$\{(M_{j_k}, g_{j_k}(t), \o_{j_k}, \vp_{j_k}(t), J_{j_k}(t), p_{j_k})\}$ which converge to $(M, g(t), \o, \vp(t), J(t), p)$.
\end{theorem}

Finally, we apply Theorem \ref{CompactnessTypeIIAFlow} to the study of singularities of the Type IIA flow. On this path, we have a new phenomena, as the formation of singularities is controlled by  $\log|\vp|_{g_\vp} +|Rm|_{g_\vp}$. Since this quantity is not homogeneous under scaling of metric, we need to distinguish two possibilities: in one case the norm $|\vp|_{g_\vp}$ may be bounded, while in the other it is not. 

\begin{theorem}\label{Models}
	Let $(M,g(t),\o,\vp(t),J(t))$ be a solution of the Type IIA flow on a compact manifold $M$, and assume $0<T\leq +\i$ be its maximum time of existence. Assume that the solution admits an injectivity radius estimate, that is $\mathrm{inj}(M,g(t))\geq{c \over |u| + |Rm|}$ for some $c>0$. Then for suitable times $t_j\to T$ and the rescalings $g_j(t)=C_jg(t_j+C_j^{-1}t)$, the sequence $g_j(t)$ admits a convergent subsequence to a solution of the Type IIB flow with a corresponding function $f(t) = \sup_M (|u| + |Rm|)$ satisfying the following properties:
	
	\smallskip
	
	{\rm [Type I]}  For the Type I the solution exists on the interval $(-\infty, c)$ and $f(t) \leq {C \over c - t}$;
	
	{\rm [Type IIa]} For the Type IIa the solution exists for all time and $f(t) \leq C$;
	
	{\rm [Type IIb]} For the Type IIb the solution exists for all time and $f(t) \leq C$ (the same as the previous case);
	
	{\rm [Type III]} For the Type III singular solution the singular model solution exists on the interval $(-C, \infty)$ and $f(t) \leq {C \over c + t}$.
	
	{\rm [Type IVa]} $\sup_{[0;T]}(T-t)|Rm| = C < \infty$, $|\vp|$ is bounded and $T<\infty$
	
	{\rm {[Type IVb]}} $\sup_{[0;T]}(T-t)|Rm| = \infty$, $|\vp|$ is bounded and $T<\infty$
	
	{\rm {[Type IVc]}} $\sup_{[0;T]}t|Rm| = \infty$,$|\vp|$ is bounded and $T=\infty$
	
	{\rm {[Type IVd]}} $\sup_{[0;T]}t|Rm| = C < \infty$, $|\vp|$ is bounded, and $T = \infty$
	
\end{theorem}

The paper is organized as follows: in Section 2, we describe the basic properties of the Type IIA geometry. In Section 3. we formulate important facts about the Type IIA flow. Sections 4 and 5 contain the proofs of Theorem \ref{CompactnessTypeIIA} and Theorem \ref{CompactnessTypeIIAFlow} respectively. Finally, Section 6 contains proof of Theorem \ref{Models}.

{\bf Acknowledgments:} I am would like to thank my advisor Professor Duong H. Phong for suggesting this problem and his constant attention, help and patience while the paper was in preparation. I also want to thank  Shuang Liang for helpful discussions, and Vladislav Cherepanov and Yulia Gorginyan for reading the text and pointing out typos.

\section{Preliminaries.}

\subsection{The Type IIA geometry.}\label{Prelim}

In this section, we recall all necessary definitions and results about the Type IIA geometry. We are following the exposition in \cite{FPPZb}.

Let $(M, \o, \vp)$ be a 6-dimensional symplectic manifold with a symplectic form $\o$ and a closed 3-form $\vp$. Assume that $\o \wedge \vp =0$. Then, as shown in \cite{FPPZb, Hit00}, one can define the map $K_\vp: TM \rightarrow \Lambda^5(M) \cong \Lambda^6(M) \otimes TM$ by the formula:

\bea
K_\vp(v) = - \iota_v \vp \wedge \vp.
\eea

We define the following "trace" of the square of $K_\vp$:

\bea
\l_\vp = {1 \over 6}\tr(K_\vp^2) \in (\Lambda^6(T^*M))^{\otimes 2}.
\eea

When $\l_\vp < 0$ we will say that $\vp$ is nondegenerate. In this case, we can extract a positive square root $\sqrt{-\l_\vp} \in \Lambda^6(M)$. Now we define the following operator:

\bea
J_\vp := {K_\vp \over \sqrt{-\l_\vp}}.
\eea
In \cite{Hit00}, Hitchin shows that $J_\vp$ is an almost complex structure. It is showed in \cite[Lemma 3]{FPPZb} that $J_\vp$ is compatible with respect to $\o$ if and only if $\o \wedge \vp =0$. This means that in this case the symmetric bilinear form $g_\vp(X,Y):=\o(X,J_\vp Y)$ is positive definite and $\o$ is invariant under $J_\vp$.

We can summarize it all in the following definition.

\begin{definition}\label{TypeIIAgeometry}
	Let $(M, \o, \vp, J_\vp, g_\vp)$, where $(M,\o)$ is a 6-dimensional symplectic manifold, $\vp$ is a closed 3-form, which satisfies $\o \wedge \vp = 0$ and $(J_\vp, g_\vp)$ are defined above. We will call such structure the Type IIA structure (or the Type IIA geometry).
\end{definition}

It is shown in \cite[Lemma 4]{FPPZb} that the following identity holds for $\l_\vp$:

\bea
\sqrt{-\l_\vp} = {1 \over 2}|\vp|^2{\o^3 \over 3!}.
\eea

It is easy to see that $d\mu_{g_\vp} = \frac{\o^3}{3!}$, where $d\mu_{g_\vp}$ is the volume form of $g_\vp$. Furthermore, we can give the explicit formula for $g_\vp$ in coordinates \cite[Lemma 5]{FPPZb}:

\bea\label{Metric}
(g_\vp)_{jk} = 12|\vp|^2\frac{\iota_{\p_j}\vp \wedge \iota_{\p_k}\vp \wedge \o}{\o^3}= - |\vp|^2\vp_{jab}\vp_{kcd}\o^{ac}\o^{bd}.
\eea 
Here and henceforth we will write $g_{jk}$ for $(g_\vp)_{jk}$.

\subsection{The Gauduchon line of connections and the holonomy of conformally rescaled metric.}

The formula \ref{Metric} implies that $\tilde{g}_{jk} = -\vp_{jab}\vp_{kcd}\o^{ac}\o^{bd}$ is also a Riemannian metric compatible with $J_\vp$. The metric $\tilde{g}_\vp = |\vp|^2g_\vp$ has rather special properties, which we are going to explain here. In order to do that, we need an additional ingredient: the Gauduchon line of connections in almost Hermitian geometry. This family of connections was introduced in \cite{Gau}, but our exposition and notations are close to \cite{FPPZb}.

Let $\tilde{\na}$ be the Levi-Civita connection for the metric $\tilde{g}_\vp$. Consider a $TM$-valued 2-form $(d^c\tilde{\o})^k_{jp} = \tilde{g}^{ks}(d^c\tilde{\o})_{sjp}$, where $\tilde{\o}(X,Y) = \tilde{g}(J_\vp X,Y)$. We can decompose ${1 \over 2}(d^c\tilde{\o})^k_{jp}$ with respect to the action of $J_\vp$ as follows:
\bea
U^k_{jp} = {1 \over 4}((d^c\tilde{\o})^k_{jp} + (d^c\tilde{\o})^k_{ab}J^a_jJ^b_p)
\\
V^k_{jp} = {1 \over 4}((d^c\tilde{\o})^k_{jp} - (d^c\tilde{\o})^k_{ab}J^a_jJ^b_p)
\\
U^k_{jp} + V^k_{jp} = {1 \over 2}(d^c\tilde{\o})^k_{jp}.
\eea
This is the standard decomposition of 2-forms on an almost complex manifold. The form $U$ is the $(1,1)$-part of ${1 \over 2}(d^c\tilde{\o})^k_{jp}$, and the form $V$ is $(2,0) + (0,2)$-part of ${1 \over 2}(d^c\tilde{\o})^k_{jp}$.  
By using this decomposition, we can define the family of unitary connections $D^t$ (which is called the Gauduchon line of connections -- see  \cite{Gau} and \cite{FPPZb}) by the following formula:

\bea
D^t_jX^k = \tilde{\na}_j X^k - N^k_{jp}X^p - U^k_{jp}X^p - tV^k_{jp}X^p,
\eea
where $N$ is the Nijenhuis tensor of $J_\vp$. The connection $D^0 = D^t \vline_{t=0}$ is called the projected Levi-Civita connection and we denote it simply by $D$. The projected Levi-Civita connection has a number of special properties in the Type IIA geometry, as it is shown in the following theorem (see \cite[Theorem 3]{FPPZb}):

\begin{theorem}\label{TypeIIA_Holonomy}
	Let $(M, \o, \vp)$ be a Type IIA geometry, and $g_\vp, \tilde{g}_\vp$ are the metrics we defined above.
	Let $D$ and  ̃$\tilde{D}$ be the projected Levi-Civita connections of $g_\vp$ and  $\tilde{g}_\vp$ respectively, $\Omega =
	\vp + iJ_\vp \vp$, and $|\Omega|_{\tilde{g}_\vp}$ the norm of $\Omega$ with respect to $\tilde{g}_\vp$. Then
	
	\smallskip
	
	{\rm (1)}  ̃$\tilde{D}({\Omega \over |\Omega|_{\tilde{g}_\vp}}) = 0$. Hence, $\tilde{D}$ has holonomy in SU(3);
	
	{\rm (2)} $D^{0,1}\Omega = 0$, so $\Omega$ is formally holomorphic;
	
	{\rm (3)} The Nijenhuis tensor of $J_\vp$ has only 6 independent components.
\end{theorem}

The theorem above shows that the Type IIA geometry is an apt analog of the Calabi-Yau geometry in the almost complex setting. Moreover, it shows the importance of the metric $\tilde{g}_\vp$. Also, it is shown in \cite{FPPZb} that the metric $\tilde{g}_\vp$ defines the norm $|\vp|^2_{g_\vp}$, and hence, metric $g_\vp$ itself. However, for our purposes it is irrelevant which metric we choose, so throughout this paper, we mainly use $g_\vp$. 

\subsection{Convergence in the sense of Gromov-Cheeger-Hamilton.}

In this section, we give a definition of convergence of Type IIA manifolds. Our definition closely follows the definitions from \cite{H95} and \cite{NK}.

\begin{definition}\label{GCH_Convergence}
	Let $\{(M_j, g_j, p_j)\}$ be a sequence of complete Riemannian manifolds, where $p_j \in M_j$ is a point.
	We say that this sequence converges to a complete pointed manifold $(M, g, p)$ in $C^k$-sense if there is a sequence of compact sets $\{U_j\}$ in $M_\infty$ with the following properties:
	
	{\rm (1)} The sequence $\{U_j\}$ is an exhaustion of $M_\infty$;
	
	{\rm (2)} There are diffeomorphisms $\Phi_j: U_j \rightarrow V_j \subset M_j$ such that $\Phi_j(p) = p_j$ and the sequence of metrics $\{\Phi_j^* g_j\}$ converges to $g_\infty$ uniformly in $C^k$-topology on compact subsets in $U_j$. 
	
\end{definition}

We also recall the theorem of Hamilton from \cite{H95}:

\begin{theorem}\label{HamiltonCompactness}
	Let $\{M_j, g_j,p_j\}$ be a sequence of complete pointed Riemannian manifolds and $\nabla_j$ be the Levi-Civita connection on $(M_j,g_j)$. Assume that the following holds uniformly in $j$:
	
	{\rm (1)} There are constants $C_k>0$, such that $|\nabla^k_jRm_j|\leq C_k$;
	
	{\rm (2)} There is a positive uniform lower bound $\iota_0$ for the injectivity radii of $M_j$ at $p_j$, $\mathrm{inj}(M_j,p_j)\geq \iota_0$.
	
	Then there is a subsequence $\{M_{j_k}, g_{j_k},p_{j_k}\}$, which converges to a complete pointed Riemannian manifold $(M_\infty, g_\infty, p_\infty)$.
\end{theorem}

This theorem is an essential ingredient in the proof of Theorem \ref{CompactnessTypeIIA}.

\section{The Type IIA flow.}

The Type IIA flow was introduced by Fei, Phong, Picard, and Zhang in \cite{FPPZb}. It is the flow of 3-forms on a 6-dimensional symplectic manifold with a Type IIA structure, which can be written as follows:

\bea\label{TypeIIAforms}
\p_t \vp = d\Lambda d(|\vp|^2\star\vp).
\eea

Here $\star$ is the Hodge star operator with respect to the metric $g_\vp$. By \cite[Theorem 1]{FPPZb} we can rewrite this flow as a Laplacian flow
\bea
\p_t \vp = -dd^\dagger(|\vp|^2\vp) + 2d(N^\dagger \cdot \vp),
\eea
where $N$ is the Nijenhuis tensor of $J_\vp$ and $(N^\dagger \cdot \vp)_{kj} = N^{a ~ b}_{~j}\vp_{akb} - N^{a ~ b}_{~k}\vp_{ajb}$. Here $N^{a ~ b}_{~j} = N^a_{~js}g^{sb}$. This formulation is similar to the Bryant's Laplacian $G_2$-flow \cite{B, BX}.

We also can rewrite the Type IIA flow in terms of metric $g_\vp$ and $u = \log|\vp|$ (see either \cite[Theorem 4]{FPPZb}, or \cite{FPPZc} for the derivation):
\bea\label{TheFlowInComponents}
\p_t g_{ij} = E_{ij} := e^{u}(-2R_{ij} + 2\na_i\na_j u -4(N^2_-)_{ij} + u_i u_j + u_{Ji}u_{Jj} + 4u_s (N^{~s}_{i~j} + N^{~s}_{j~i}))
\\
(\p_t - e^{u}\D)u = (2|\na u|^2 + |N|^2)
\eea

In \cite[Proposition 7 and Theorem 7]{FPPZb} the authors proved the Shi-Type estimates for the Type IIA flow. 
\begin{theorem}\label{Shi-typeTypeIIA}
	Let $(g(t), u(t))$ evolve along the Type IIA flow. Denote by $Rm$ the curvature tensor of the Levi-Civita connection with respect to $g_\vp$, and $u=\log|\vp|$ as usual. Assume that $|u|+ |Rm| \leq C$ for some uniform constant $C > 0$ on the interval $[0;T]$. Then the following holds
	
	\smallskip
	
	{\rm(1)} There is a constant $C_0$ such that 
	\bea
	\sup_{M \times [0;T]}(|Rm| + |N|^2 +  |\na N| + |\na^2u| + |\na u| + |u|) \leq C_0
	\eea
	
	{\rm(2)} There are constants $C_k$ for $k \in \N$, such that
	\bea
	\sup_{M \times [0;T]}(|\na^k Rm| + |\na^{k+1}N| + |\na^{k+2}u|) \leq C_k
	\eea
	
	{\rm(3)} Assume that the bound $|u|+ |Rm| \leq C_0$ holds on $[0;T)$. Then $|\na^k \vp| \leq C(C_0, k, T, \vp(0))$ and the flow can be continued to $[0;T+\ve)$ for some $\ve > 0$.
\end{theorem}

The third statement of Theorem \ref{Shi-typeTypeIIA} justifies the following definition.

\begin{definition}\label{Singularities_defn}
	Consider the solution $(M, g(t), \o, \vp(t), J(t))$ for the Type IIA flow on a finite time interval $[0;T)$ for $T < +\i$. We say that the solution of the Type IIA flow becomes singular at time $T$ if $ \sup_{M \times [0;T)}(|u| + |Rm|) = + \i$. 
\end{definition}

\begin{remark}
	It is worthy to mention that in \cite{_FP21_} the authors showed that the Type IIA flow is very close to the Hitchin flow (see \cite{Hit00,_Fei15_} for the details). However, the Hitchin flow is more degenerate than the Type IIA flow, and it is hard to prove the short-time existence of the solution for the Hitchin flow. Also, the evolution equation for $\log|\vp|^2$ implies that the Hitchin functional $H(\vp) = \int \vp \wedge \hat{\vp}$ is monotone increasing along the Type IIA flow. 
\end{remark}

\section{Proof of Theorem 1.}

\medskip
Here we give the proof of Theorem \ref{CompactnessTypeIIA}.

\medskip
\textbf{Step 1: Compactness for the underlying Riemannian manifolds.}

This is the content of Theorem \ref{HamiltonCompactness}. Hence, there exists a complete, pointed Riemannian manifold $(M, g, p)$, such that $(M_j, g_j,p_j)$ converge to it in the sense of Definition \ref{GCH_Convergence}.

\medskip
\textbf{Step 2: Convergence of complex structures.}

By the formula $\na^{g_j}\o_j = - 2N_j$ (this follows from the formula for the Gauduchon connection for $g_j$) we see that for any $k \geq 1$ the norms of covariant derivatives $|(\na^{g_j})^k\o_j|_{g_j}$ are bounded. We will argue like in \cite{LotayWei15} and \cite{NK}.
Let $\Phi_{j,k}$ be restriction of $\Phi_{j+k}$ on $U_j$. Denote $g_{j,k} := \Phi_{j,k}^*g_{j+k}$. The metrics $g_{j,k}$ converge to $g_{j, \i}$ as $k \rightarrow \i$. It is easy to see that $g_{j, \i}$ is the restriction of $g_\i$ on $U_j$. Similarly, the notation $\o_{i,j}$ means $\Phi^*_{i,j}\o_j$ and $\vp_{i,j}$ stands for $\Phi^*_{i,j}\vp_j$.

Consider the quantity $A_{i,j} = \na^g - \na^{g_{i,j}}$. It is bounded together with all covariant derivatives with respect to $\na^g$, because $g_{i,j}$ converges to $g$ in $C^\i$. Using this, together with the fact that $g_{i,j}$ is uniformly equivalent to $g$ for big $j$, we obtain that there are constants $C'_k$ such that $|(\na^g)^k\o_{i,j}|_g \leq C'_k$. By the Arzela-Ascoli theorem, after passing to a subsequence, there is a form $\o_i$, such that the subsequence converges in $C^\i$ to $\o_i$. The form $\o_i$ satisfies the identity $\o_i^3 = 6d\mu_g$, hence it is nondegenerate. Moreover, there is a global form $\o$ on $M$, such that $\o|_{U_i} = \o_i$. 

Finally, since $J_j = g^{-1}_j \o_j$, we see that $J_{i,j} = g^{-1}_{i,j}\o_{i,j}$ converge to the operator $J_{i,\i} = g^{-1}\o_i$, which satisfies $J^2_{i, \i} = -\mathrm{Id}$. Moreover, all $J_{i,\i}$ are gluing together to a globally-defined almost complex structure $J = g^{-1}\o$.

\medskip
\textbf{Step 3: Convergence of 3-forms.}

The argument that proves the convergence of $\vp_{i,j}$ to $\vp_{i, \i}$ and the existence of a globally defined limit form $\vp$ carries over verbatim. We only need to prove the nondegeneracy of the form $\vp$.

We use the following identities in Type IIA geometry (see Section \ref{Prelim}): 
\bea\label{ThreeForm}
g(\xi_1,\xi_2) {\o^3 \over 3!} = \iota_{\xi_1}\vp \wedge \iota_{\xi_2}\vp \wedge \o
\\
K_\vp\xi = - \iota_{\xi}\vp \wedge \vp
\\
J\xi = {K_\vp\xi \over \sqrt{-\tr(K^2)}}
\\
\sqrt{-\tr(K^2)} = {1 \over 2} |\vp|^2 {\o^3 \over 3!}
\\
\Lambda_{\o} \vp = \o \wedge \vp = 0
\eea  

These identities hold for any tuple $(g_{i,j},\o_{i,j},J_{i,j},\vp_{i,j})$ and for any vector fields $\xi, \xi_1$ and $\xi_2$.

Again, by the Arzela-Ascoli theorem, there are limiting 3-forms $\vp_i$ on $U_i$, and these forms define a global form $\vp$ on $M$. Next, the positive lower bound for $|\vp_{i,j}|$ implies that $|\vp| > 0$ pointwise, and ${1 \over 2} |\vp|^2 {\o^3 \over 3!} > 0$. This implies that all the identities \ref{ThreeForm} holds for $(g,\o,J,\vp)$. Hence, the form $\vp$ is a nondegenerate 3-form, and the tuple $(g,\o,J,\vp)$ defines the Type IIA geometry on $M$. Q.E.D.

\section{Compactness of solutions of the Type IIA flow.}

\subsection{Estimates for mixed derivatives.}

In this section we are going to prove the estimates for $\p_t^m\na^kg_j$ and $\p_t\na^k \vp_j$. These estimates are the same as in the case of the Ricci flow, the Laplacian $G_2$ flow, and the Type IIB flow.

Next, we prove that metrics are equivalent under the assumptions of Theorem \ref{CompactnessTypeIIAFlow}:
\bea
C_1 := C|t_1 - t_0| \geq \int_{t_0}^{t_1}|\p_t \log g(t)(V,V)|dt \geq |\log{g(t_1)(V,V) \over g(t_0)(V,V)}|,
\eea
where $C$ depends on the bounds for $Rm, N$ and $u$. This implies that norms defined by $g(t_0)$ and $g_(t_1)$ are equivalent. 

Consider $\O \times [a;b] \subset M \times (\a, \b)$ such that $\O \subset U_j$ for j large enough.  By abusing the notation, we will denote $\Phi_j^*\vp_j(t)$ and $\Phi_j^*g_j(t)$ on $\O$ by $\vp_j$ and $g_j$ respectively. To stress the difference between $g_j$ and $g$ on $M$ we will denote $g$ by $\hat{g}$ and $\vp, \o$ by $\hat{\vp}$ and $\hat{\o}$ respectively.

Now we prove the estimates for $|\p^l_t\hat{\na}^kg_j|_{\hat{g}}$ and $|\p^l_t\hat{\na}^k\vp_j|_{\hat{g}}$. To begin with, consider $|\p_t\hat{\na}g_j|_{\hat{g}}$. Let $A_j := \hat{\na} - \na^{g_j}$ and $E_j = \p_t g_j(t)$. We have:
\bea
|\p_t\hat{\na}g_j|_{\hat{g}} \leq |\hat{\na}\p_t g_j|_{\hat{g}} \leq |A_j|_{\hat{g}}|E_j|_{\hat{g}} + |\na^{g_j}E_j|_{\hat{g}} \leq C(1 + |\hat{\na}g_j|_{\hat{g}})
\eea
Here we can estimate the norm $|E_j|_{\hat{g}}$ because the metrics $g_j$ and $\hat{g}$ are equivalent. Using the following formula, we can proceed by induction:
\bea
\p_t\hat{\na}^kg_j = \sum_{j_1 + \dots + j_p = k} \na^{j_1} E_j * \hat{\na}^{j_2} g_j * \dots * \hat{\na}^{j_p} g_j
\eea
which can be proved in a similar manner as in the case of the Ricci flow (see \cite[Lemma 8.6]{AndrewsHopper}). Now the induction by $k$ shows that for each $k$ there is a constant $C$, independent of $j$, such that
\bea
|\p_t\hat{\na}^kg_j|_{\hat{g}} \leq C(1 + |\hat{\na}^kg_j|_{\hat{g}}),
\eea
so this implies that $|\p_t\hat{\na}^kg_j|_{\hat{g}}$ is bounded.
Finally, note that $\p_t^lg_j = \p_t^{l-1} E_j$, and formula \ref{TheFlowInComponents} together with evolution formulas in \cite[Proposition 7 and Theorem 7]{FPPZb} imply that the covariant derivatives of $\p_t^{l-1} E_j$ are uniformly bounded.

The similar arguments holds for $|\p^l_t\hat{\na}^k\vp_j|_{\hat{g}}$. Consider the case $l=1$. We have:
\bea
|\p_t\hat{\na}^k\vp_j|_{\hat{g}} = |\hat{\na}^k\p_t\vp_j|_{\hat{g}} = |\hat{\na}^k(dd^\dagger(|\vp_j|_{j}^2\vp_j) + 2d(N_j^\dagger \cdot \vp_j))|_{\hat{g}} \leq C(|\hat{\na}^{k+2}\vp_j|_{\hat{g}} + |\hat{\na}^{k+1}\vp_j|_{\hat{g}}).
\eea
By the equivalence of metrics and Theorem \ref{Shi-typeTypeIIA}, the RHS is bounded. For $l>1$ we have
\bea
|\p^l_t\hat{\na}^k\vp_j|_{\hat{g}} \leq |\p^{l-1}_t\hat{\na}^kdd^\dagger(|\vp_j|_{j}^2\vp_j)|_{\hat{g}} + 2|\p^{l-1}_t\hat{\na}^k(d(N_{\vp_j}^\dagger \cdot \vp_j))|_{\hat{g}} \leq
\\
C(|\p^{l-1}_t\hat{\na}^{k+1}(|\vp_j|_{j}^2\vp_j)|_{\hat{g}} + |\p^{l-1}_t\Sigma_{j=0}^{k+1} \hat{\na}^jN^\dagger \hat{\na}^{k-j+1}\vp_j|_{\hat{g}})
\eea 
The second summand can be estimated in terms of $|\p^m_t\hat{\na}^p\vp_j|_{\hat{g}}$ and $|\p^m_t\hat{\na}^p N_j|_{\hat{g}}$ for $m \leq l-1$. However, the Nijenhuis tensor $N_j$ can be expressed via $J_j$ and $\hat{\na}J_j$, while the time derivatives of $J_j$ can be expressed in terms of time derivatives of $g_j$. More precisely, $\p^l_t J^a_b = \o^{as}\p^{l-1}_t(\p_t g_{sb})$. Hence, the second summand in the formula above can be bounded by a finite expression, involving $|\p^m_t\hat{\na}^p\vp_j|_{\hat{g}}$ and $|\p^m_t\hat{\na}^pg_j|_{\hat{g}}$.  The induction by $l$ shows that the norms $|\p^l_t\hat{\na}^k\vp_j|_{\hat{g}}$ are also bounded for each $l>1$. 

Thus we proved the following proposition

\begin{proposition}\label{MixedDerivatives}
	Under the assumptions of Theorem \ref{CompactnessTypeIIAFlow} there are constants $C_{l,k}$, such that $|\p^l_t\hat{\na}^kg_j|_{\hat{g}} \leq C_{l,k}$ and $|\p^l_t\hat{\na}^k\vp_j|_{\hat{g}} \leq C_{l,k}$ on $\O \times [a;b] \subset M \times (\a, \b)$.
\end{proposition}

\subsection{Proof of Theorem \ref{CompactnessTypeIIAFlow}.}

First of all, Theorem \ref{CompactnessTypeIIA} guarantee that, after taking a subsequence, there is a limit $(M, g, \o, \vp, J, p)$ of $(M_j, g_j(0), \o_j, \vp_j(0), J_j, p_j)$. 

Next, by Proposition \ref{MixedDerivatives} and the Arzela-Ascoli theorem (see \cite[Corollary 9.14]{AndrewsHopper}), the sequence $\{(g_j(t), \vp_j(t), J_j(t))\}$ converge (again after taking a subsequence) to the triple $(g(t), \vp(t), J(t))$ on $\O \times [a;b] \subset M \times (\a, \b)$ in $C^\i$. The $C^\i$ convergence implies that the couple $(g(t), \vp(t))$ satisfies the Type IIA flow equation. Now the standard diagonalization argument implies that there exists $(M, g(t), \vp(t), J(t), \o)$ for all times $t \in (\a, \b)$, where $\vp(t), \g(t)$ are the solutions to the Type IIA flow. Q.E.D.

\section{The singularities of the Type IIA flow.}

In this section, we describe the singularities for the Type IIA flow. Recall that from Definition \ref{Singularities_defn} singularities are characterized by the condition $\sup_{M \times [0;T)}(|u| + |Rm|) = +\i$.
From this, we can see a significant difference between the Type IIA flow and many other geometric flows such as the Ricci flow, the Type IIB flow, and the Laplacian $G_2$-flow. The singularities of all these flows are characterized by quantities, which are {\it homogeneous} under the rescaling of metric. However, in the case of the Type IIA flow, we have a different situation. If we multiply metric $g$ by a positive constant, then the norm of the Riemann curvature tensor and the function $|u| = |\log |\vp|^2|$ behave differently under this operation. Hence, the usual approach to the singularities doesn't work here and we need to consider another quantity, which behaves nicely under the rescaling.

Consider the function $f(t):=\sup_M(|\vp|^{1 \over 3} + |Rm|)$. This quantity has the desired homogeneous properties but it is now clear whenever $f(t)$ has the same rate of growth as $|u| + |Rm|$. There are two possible cases:
\medskip

{\rm 1.} The norm $|\vp|$ is bounded along the flow;

{\rm 2.} The norm $|\vp|$ tends to $+\i$ as $t \rightarrow T$.
\medskip

The latter case is simpler: if $|\vp| \nearrow +\i$, we can find a sequence of times $t_j \rightarrow T$, such that $(|u| + |Rm|)(t_j) \leq f(t_j)$. The former case is a bit more delicate: if $|\vp|$ remains small enough along the flow, we cannot say that $|\vp|^{1 \over 3} \geq \log|\vp|^2$. However, if $|\vp|$ is bounded, the condition $\sup_{M \times [0;T)}(|u| + |Rm|) = + \i$ is equivalent to the condition $\sup_{M \times [0;T)} |Rm| = + \i$.

With all these considerations we can give the following definition.

\begin{definition}\label{SingularityTypes}
	Suppose $(M, g(t), \o, \vp(t), J)$ be a solution of the Type IIA flow. Let $f(t)$ be as above, and $T$ is the maximal time for which flow exists. We define the following type of singularities
	
	{\rm [Type I]} $\sup_{[0;T]}(T-t)f(t) = C < \infty$, and $T<\infty$	
	
	{\rm [Type IIa]} $\sup_{[0;T]}(T-t)f(t) = \infty$, and $T<\infty$
	
	{\rm [Type IIb] }$\sup_{[0;T]}tf(t) = \infty$, $T=\infty$
	
	{\rm [Type III]} $\sup_{[0;T]}tf(t) = C < \infty$, $T = \infty$	
	
	{\rm [Type IVa]} $\sup_{[0;T]}(T-t)|Rm| = C < \infty$, $|\vp|$ is bounded and $T<\infty$
	
	{\rm {[Type IVb]}} $\sup_{[0;T]}(T-t)|Rm| = \infty$, $|\vp|$ is bounded and $T<\infty$
	
	{\rm {[Type IVc]}} $\sup_{[0;T]}t|Rm| = \infty$,$|\vp|$ is bounded and $T=\infty$
	
	{\rm {[Type IVd]}} $\sup_{[0;T]}t|Rm| = C < \infty$, $|\vp|$ is bounded, and $T = \infty$
	
\end{definition}

\begin{remark}
	Here we have a problem with terminology: we study the Type IIA flow, but the singularity types are named like Type I, Type II, etc. Our terminology follows closely the terminology from \cite{Hamilton1}, and we hope it is not misleading. 
\end{remark}


With all these preparations, we can prove Theorem \ref{Models}.

\medskip
\textbf{Proof of Theorem \ref{Models}}

The proof is essentially the same as the proof of the similar theorem for the Ricci flow and \cite[Theorem 3]{NK}. 

Consider the family of metrics $g_j(t) := C_jg(t_j + {t \over C_j})$, where we determine the constants $C_j$ and $t_j$ later. In addition, we rescale the symplectic form $\o$ by $\o_j = C_j\o$. Also, we need to rescale the form $\vp$ as following: $\vp_j(t) = \vp(t_j + {t \over C_j})$.  

\medskip

\noindent{\em The case of Type I singularity}

\medskip

Let $C = \sup_{[0;T]}(T-t)f(t)$ and pick a sequence $(x_j,t_j) \in M \times [0;T)$, and assume $C_j:=F(x_j,t_j)$. Note that $(T - t_j)C_j \geq c >0$ for some $c$. Now we consider rescaling of the metric and forms as described above. We have
\bea
|\vp_j|^{1 \over 3}_{g_j} + |Rm_j|_{g_j} \leq C_j^{-1}{C \over T - (t_j + C_j^{-1}t)} \leq {C \over c -t}
\eea

Now we applying the compactness theorem and extract a limit of the sequence \\ $(M, g_j(t), \o_j, \vp_j(t), J_j)$. 

\medskip

\noindent{\em The case of Type IIa singularity}

\medskip

The proof is again an adaptation of the proof in the case of the Ricci flow. Pick $T_j \nearrow T$ and points $(x_j, t_j)$ such that $(T_j - t_j)F(x_j,t_j) = \max_{M \times [0;T_j]}F(x,t)$. Taking $C_j$ as before, we have 

\bea
|\vp_j|^{1 \over 3}_{g_j} + |Rm_j|_{g_j} = {1 \over C_j}(|\vp| + |Rm|) \leq {C_j (T_j - t_j) \over C_j (T_j - t_j - {t \over C_j})} = {1 \over 1 - {t \over C_j(T_j - t_j)}}\leq 1
\eea

\medskip

\noindent{\em The case of Type IIb singularities.}

\medskip

We again pick up the sequence of points $(x_j,t_j)$, such that $t_j(T_j - t_j)F(x_j,t_j) = \max_{M \times [0;T_j]}t(T_j - t)F(x,t)$.
By rescaling, we obtain
\bea
|\vp_j|^{1 \over 3}_{g_j} + |Rm_j|_{g_j} = {1 \over C_j} (|\vp|^{1 \over 3} + |Rm_j|)(t_j + {t \over C_j}) \leq {C_jt_j(T_j - t_j) \over C_j(t_j + {t \over C_j})(T_j - t_j - {t \over C_j})} = 
\\
={1 \over (1 + {t \over C_jt_j})(1 - {t \over C_j(T_j - t_j)})} = {C_jt_j \over C_jt_j - t}{C_j(T_j - t_j) \over C_j(T_j - t_j) - t} \rightarrow 1
\eea

Applying the compactness theorem, we obtain the desired statement.

\medskip

\noindent{\em The case of Type III singularity}

\medskip

As in the case of Type I singularities, let $C = \sup_{M \times [0;+ \i)}tf(t)$, and $(x_j, t_j) \in M \times [0;\i)$ as before. Then

\bea
|\vp_j|^{1 \over 3}_{g_j} + |Rm_j|_{g_j} \leq {C \over C_j (t_j + {t \over C_j})} \leq {C \over c + t},
\eea
since $C_jt_j \geq c > 0$. Hence, after passing to a subsequence and applying Theorem \ref{CompactnessTypeIIAFlow}, the limit satisfies the inequality $|\vp|^{1 \over 3} + |Rm| \leq {C \over c + t}$.
\medskip

\noindent{\em The case of the Type IV singularities.}
\medskip

Let and assume that $g_j(t)$ and $\o_j$ are the same as before. We also rescale $\vp_j =  C^{3 \over 2}_j\vp(t_j - {t \over C_j})$. Then $|\vp_j|_{g_j}(t)$ remains bounded, so we can apply the compactness theorem and extract the limit. The classification part carries out verbatim as in \cite{H95}. Q.E.D. 

Finally, we consider a few examples from \cite{FPPZb} and  

\begin{example}[The Type IIa singularity]
	Consider the 6-dimensional solvmanifold $M$ with the complex structure, constructed by Tomassini and Vezzoni in \cite{TV}. The invariant forms satisfy the following identities:
	\bea
	de^1 = -\l e^{15}, ~de^2= \l e^{25}, ~de^3 = -\l e^{36}
	\\
	de^4 = \ e^{46}, ~de^5 = 0, ~de^6 = 0.
	\eea
	and the symplectic form is given by $\o = e^{12}+ e^{34} + e^{56}$. Here $e^{ij} = e^i \wedge e^j$, and $\l = \log {3 + \sqrt{5} \over 2}$. 
	
	The nondegenerate 3-form can be defined as follows:
	\bea
	\vp = \a(e^{135} + e^{136}) + \b(e^{145} - e^{146}) + \g(e^{235} - e^{236}) - \d(e^{245} - e^{246}),
	\eea
	and $\a,\b,\g,\d$ depends on $t$. The nongedeneracy condition is $|\vp|^4 = 64\a\b\g\d > 0$. 
	
	In \cite{FPPZb} the authors showed that the Type IIA flow preserves this ansatz and exists for a finite time $T = {1 \over 32 \l^2}{\log(\b_0\g_o) - \log(\a_0\d_0)\over\b_0\g_o - \a_0\d_0}$, which is defined by the initial data. Moreover, the norm of $\vp$ tends to infinity. However, the metric $g_\vp$, the complex structure $J_\vp$ and the Nijenhuis tensor $N$ extends smoothly to $t=T$, hence the curvature is bounded as $t \rightarrow T$. Using the calculations from \cite{FPPZb}, we can directly compute $|u|$:
	\bea
	|\log |\vp|^2| = |\log(8\sqrt{\a\b\g\d})| = \left|-16\l^2t + \log {8|\b_0\g_0 - \a_0\d_0|\sqrt{\a_0\b_0\g_0\d_0}\over|\b_0\g_0e^{\a_0\d_0} - \a_0\d_0e^{\g_0\b_0}|}\right|.
	\eea 
	
	Hence, ${|u|\over T-t}$ tends to infinity unless $T = 16\l^2\log {8|\b_0\g_0 - \a_0\d_0|\sqrt{\a_0\b_0\g_0\d_0}\over|\b_0\g_0e^{\a_0\d_0} - \a_0\d_0e^{\g_0\b_0}|}$. In the latter case, ${|u| + |Rm|\over T-t}$ is bounded, so this is the Type I singularity.
\end{example}

\begin{example}[The Type IIb singularity]
	Consider the nilmanifold $M$ from \cite[Example 5.2]{BT} (see also \cite[Section 9.3.2]{FPPZb}). The Lie algebra of the nilpotent Lie group can be characterized by invariant 1-forms:
	\bea
	de^1 = de^2 = de^3 = de^5 = 0
	\\
	de^4 = e^{15}, de^6 = e^{13}.
	\eea
	The symplectic form is standard: $\o = e^{12} + e^{34} + e^{56}$.
	The family of nondegenerate primitive 3-forms can be defined as follows:
	\bea
	\vp_{a,b} = (1+a)e^{135} - e^{146} - e^{245} - e^{236} + b(e^{134} - e^{156}),
	\eea
	where $a$ and $b$ are functions of $t$. From \cite{FPPZb} we know, that the positivity condition means $|\vp_{a,b}|^4 = 16(1 + a -b^2) > 0$, and in this case the flow exists for all time $t>0$. Moreover, it was shown that $a(t) = a_0 + 8t$ and $b=b_0$, where $a_0, b_0$ are initial data.
	
	Since $|u| = |\log4(1 + 8t + a_0 - b^2_0)|$, we see that $t(|u| + |Rm|) \rightarrow +\i$ as $t \rightarrow +\i$. Hence, the solution forms Type IIb singularity at $t=+\i$.
\end{example}

\bigskip

\noindent Department of Mathematics, Columbia University, New York, NY 10027 USA

\noindent nklemyatin@math.columbia.edu

\end{document}